\documentclass[10pt]{amsart}
\usepackage[english]{babel}
\usepackage{amssymb,amsmath,amsthm,amscd,mathrsfs}
\usepackage{color}   
\newcommand{\ec}{\color{black}}%
\newfont{\rams}{msbm10 scaled\magstep1}
\newfont{\ramss}{msbm10 scaled\magstep0}
\newfont{\iams}{msbm10}
\newfont{\gotic}{eufm10 scaled\magstep1}
\newfont{\bellap}{eusm10 scaled\magstep1}

\newcommand{\im}{{\rm Im}}

\newcommand{\sss}{{\rm Spec}}

\newcommand{\zar}{{\rm Zar}}

\newcommand{\ff}{\mathscr F}
\newcommand{\uu}{\mathscr U}
\newcommand{\ms}{\mathscr}

\newcommand{\ad}{\mbox{\it\texttt{Cl}}}

\newcommand{\z}{{\ldots}}
\newcommand{\w}{{\setminus}}

\swapnumbers
\newtheoremstyle{break}
  {9pt}
  {9pt}
  {\itshape}
  {}
  {\textsc}
  {.}
  {.7em}
  {}
\newtheoremstyle{break1}
  {9pt}
  {9pt}
  {}
  {}
  {\textsc}
  {.}
  {.7em}
  {}
\theoremstyle{break}
\newtheorem{thm}{ \textsc{Theorem}}[section]

\newtheorem{cor}[thm]{ \textsc{Corollary}}
\newtheorem{lem}[thm]{ \textsc{Lemma}}
\newtheorem{prop}[thm]{ \textsc{Proposition}}
\theoremstyle{break1}

\newtheorem{oss}[thm]{ \textsc{Remark}}
\theoremstyle{remark}

\title[Ultrafilter and constructible topologies on spaces of valuations]{Ultrafilter and Constructible topologies\\ on spaces of valuation domains}
\dedicatory{\sl In memory of Nicolae Popescu}
\author{Carmelo A. Finocchiaro, Marco Fontana, and K. Alan Loper}
\address{C.F. \& M.F., Dipartimento di Matematica, Universit\`a degli Studi
``Roma Tre'', Largo San Leonardo Murialdo, 1,  00146 Rome, Italy.}
\email{carmelo@mat.uniroma3.it }
\email{fontana@mat.uniroma3.it }
\address{K.A.L., Department of Mathematics, Ohio State University, Newark, OH 43055, USA}
\email{lopera@math.ohio-state.edu}

\thanks{\it Acknowledgments. \rm
During the preparation of this paper, the first two authors
were  partially supported by  a research grant PRIN-MiUR}

\date{December 16, 2011}

\begin{document}

\maketitle

\ec
\hfill{\footnotesize \sl }%

\bigskip

\begin{abstract} Let $K$ be a field and let $A$ be a subring of $K$. We consider pro\-per\-ties and applications of a compact, Hausdorff topology
 called the ``ultrafilter topology''  defined on the space Zar$(K|A)$ of all valuation domains having $K$ as quotient field and containing $A$.  We show that the ultrafilter topology coincides with the constructible topology on the abstract Riemann-Zariski surface Zar$(K|A)$. We extend results  regarding distinguished spectral topologies on spaces of valuation domains.
\end{abstract}

\medskip

\section{Introduction}
Let $K$ be a field and let $A$ be a subring of $K$.  
 We denote  by Zar$(K|A)$ the collection of all valuation domains which have $K$ as quotient field and  have $A$ as a subring.  In case $A$ is the prime subring of $K$, then Zar$(K|A)$ includes all valuation domains with $K$ as quotient field and we denote it by simply 
Zar$(K)$.  The  first topological approach to the space Zar$(K)$ is due to Zariski who proved the quasi-compactness of this space, endowed with what is now called the Zariski topology (see \cite{za} and \cite{zasa}).
Later, it was proven, and rediscovered by several authors with a variety of different techniques, that if $K$ is the quotient field of $A$ then Zar$(K|A)$ endowed with Zariski's topology  is a spectral space in the sense of Hochster \cite{ho} (see \cite{dofefo}, \cite{dofo}, \cite{hu-kn} and the appendix of \cite{ku}).

  In Section 2, we define the Zariski topology on Zar$(K)$  and a classical refinement of it known as the constructible topology.  We also introduce the notion of an ultrafilter and point out that the current authors recently used ultrafilters to define a topology on the set Spec$(R)$ of prime ideals of a commutative ring and then prove that this ultrafilter topology is identical  with the classical constructible topology on Spec$(R)$ \cite{folo}.  In Section 3, we define a constructible topology and an ultrafilter topology on the space Zar$(K|A)$ for {\sl any} subring $A$ of $K$ and demonstrate that they are identical.   
    In a subsequent paper we will study further the ultrafilter/constructible topology on the space Zar$(K|A)$ providing some applications to  the representations of integrally closed domains as intersections of va\-lua\-tion overrings  \cite{fi-fo-lo-11} (see also \cite{Finocchiaro}).

  This paper is dedicated to the memory  of Nicolae Popescu who recently left us: his papers were an important source of inspiration (e.g., \cite{AP}, \cite{APZ} \cite{PV-1}, \cite{PV-2},  \cite{PZ},  and \cite{KMPV}).

\medskip

\section{ Notation and preliminaries}
\label{not}
If $X$ is a set, we  denote by $\boldsymbol{\mathscr{B}}(X)$ the collection of all subsets of $X$, and by $\ms B_{\mbox{\Tiny{\texttt{fin}}}}(X)$ the collection of all finite subsets of $X$. Moreover, if $\mathscr G$ is a nonempty subset of $\ms B(X)$,  we will { simply} denote by $\bigcap \mathscr G$ { the set obtained by intersection} of all {  subsets of $X$ belonging to} $\mathscr G$, i.e., $\bigcap \mathscr G := \bigcap \{ G \mid G \in  \mathscr G\}$.

 Recall that a nonempty collection $\ff$ of subsets of $X$ is said to be a \emph{filter on}  $X$ if the following conditions are satisfied:
\rm(a) $\emptyset \notin \ff$; 
\rm(b) if $F,G\in \ff$, then $F\cap G\in \ff$;
\rm(c) if $F,G\in \ms B(X)$, $F\subseteq G$, and $F\in \ff$, then $G\in \ff$.

Let $\boldsymbol{\mathcal F}(X)$ be the set of all filters on $X$, partially ordered by inclusion.  We say that a filter $\ff$ on $X$ is an  \emph{ultrafilter on} $X$ if it is a maximal element in $\boldsymbol{\mathcal F}(X)$.  In the following, we  denote the collection of all ultrafilters on a set  $X$ by $\boldsymbol{\beta}(X)$.

For each  $x\in X$,  it is immediately seen that $\beta_X^x:=\beta^x:=\{Z\in\ms B(X)\mid x\in Z\}$ is an ultrafilter on $X$, called the \emph{trivial \emph{(or} fixed \emph{or} principal\emph{)} ultrafilter of $X$ {centered on $x$}}.

In the next lemma, we collect  some basic facts of filters and ultrafilters needed in this paper.

\begin{lem}\label{siparte}
Let $X$ be a set. 
\begin{enumerate}
\item[\rm(1)] If $\ff$ is a filter on $X$, then there is an ultrafilter $\uu$ on $X$ such that $\ff\subseteq \uu$.

\item[\rm(2)] If $\mathscr G$ is a collection of subsets of $X$ with the finite intersection property, then there is a filter $\ff$ on $X$ such that $\mathscr G\subseteq \ff$.

\item[\rm(3)] Let $f:X\rightarrow Y$ be a map and $\ms U$  an ultrafilter [respectively, $\ff$ a filter] on $Y$.  If $f$ is injective and $f(X) \in \uu$ [respectively, $f(X) \in \ff$],  then 
$$\ms U^{{\!f}}:=\{f^{-1}(Z) \mid  Z \in \ms U\}
 \; \mbox{[respectively, $ \ff^{{\!f}}:=\{f^{-1}(Z) \mid Z\in \ff\}$]} 
$$
 is an ultrafilter [respectively, a filter] on $X$. In particular, if $X$ is a subset of $Y$ and $f$ is the inclusion map, then the set
$$
\uu^{{\!X}}:=\{Z\cap X \mid Z\in \uu\} \; 
\mbox{[respectively, $ \ff^{{\!X}}:=\{Z\cap X \mid Z\in \ff\}$]} 
$$
is an ultrafilter [respectively, a filter] on $X$. { Moreover,} in this case,  $ \uu^{{\!X}} \subseteq \ms U$ [respectively, $ \ff^{{\!X}}\subseteq \ms F$].
\item[\rm(4)] Let $f:X\rightarrow Y$ be a map and let $\ms U$ be an ultrafilter [respectively, $\ff$ be a filter]  on $X$, then 
$$
\begin{array} {rl}
\ms U_{_{\!f}} & \hskip -10 pt :=\{Z\in\ms B(Y) \mid f^{-1}(Z)\in \ms U\} \\
 & \hskip -10 pt \mbox{[respectively, $ \ff_{_{\!\!\!f}}:=\{Z\in \ms B(Y) \mid f^{-1}(Z) \in \ff \}$]} 
 \end{array}
$$
 is an ultrafilter [respectively,  a filter]  on $Y$. In particular, if $X$ is a subset of $Y$, $f$ is the inclusion map  and $\ms U$ is an ultrafilter [respectively, $\ff$ is a filter] on $X$, then the set
$$
\begin{array} {rl}
 {\ms U}_{_{\!Y}} & \hskip -10 pt:=\{Z\in \ms B(Y) \mid Z\cap X\in \ms U\}\\
  & \hskip -10 pt \mbox{[respectively, $ \ff_{_{\!Y}} :=\{Z\in \ms B(Y) \mid Z\cap X\in \ff \}$]} 
   \end{array}
$$
is an ultrafilter [respectively,  a filter]  on $Y$. { Moreover,} in this case,  $ \ms U \subseteq   {\ms U}_{_{\!Y}}$ [respectively, $ \ms F \subseteq {\ms F}_{_{\!Y}}$].

\item[\rm(5)] If $\ff$ is a filter on $X$, then the following conditions are equivalent.
\begin{enumerate}
\item[\rm(i)] $\ff$ is an ultrafilter.
\item[\rm(ii)] If $Y,Z\in \ms B(X)$ and $Y\cup Z\in \ff$, then either $Y\in \ff$ or $Z\in \ff$.
\item[\rm (iii)] If $Y\in \ms B(X)$, then either $Y\in \ff$ or $X\w Y\in \ff$.
\end{enumerate}
\end{enumerate}
\end{lem}
\noindent
\textsc{Proof}. (1) is proved in  \cite[Theorem 7.5]{jech}. (2) Note that the collection
$$
{ \ff({\mathscr G})}:=\{Z\in\mathscr B(X) \mid Z\supseteq \bigcap \mathscr G', \mbox{  for some } \mathscr G'\subseteq \mathscr G,\,\mathscr G' \mbox{ finite} \}
$$
is a filter on $X$ and, precisely, it is the smallest filter on $X$ containing $\mathscr G$ (see also \cite[Lemma 7.2(iii)]{jech}). 
(3) {is an easy consequence of definitions} and \cite[Exercise 7.1]{jech}. The first part of (4) { is given in} \cite[Exercise 7.5]{jech}. The second part of (4) { is a straightforward consequence of the first one}. {Finally,} (5) is proved in \cite[Lemma 7.4 and Exercise 7.3]{jech}. \hfill$\Box$

\medskip

\medskip

If $K$ is a field and $A$ is a subring of $K$, then we denote by $\zar(K|A)$ the set of all valuation rings of $K$ containing $A$, and simply by $\zar(K)$  the set $\zar(K|A_1)$  when $A = A_1$ is the fundamental subring of $K$. 

As  is well known, Zariski { \cite{za} (or, \cite[Volume II, Chapter VI, \S 1, page 110]{zasa})} introduced and studied the set ${ Z}:=\zar(K|A)$ together with a topological structure defined by taking, as a basis for the  open sets, the subsets  
$
B_{{\!F}}^{{ Z}} := \{V\in {Z}  \mid V \supseteq F\}
$, for $F$ varying in $ \ms B_{\mbox{\Tiny{\texttt{fin}}}}(K)$, i.e., if $F:=\{x_1, x_2, \dots, x_n\}$, with $x_i \in K$, then $$
B_F^{{Z}} = \zar(K|A[x_1, x_2, \dots, x_n]).$$
This topology is called the \emph{Zariski topology on}  $Z=\zar(K|A)$ and $Z$, equipped with this topology, also denoted  later by ${{ Z}}^{\mbox{\tiny{\texttt{zar}}}}$, is usually called the \emph{(abstract) Zariski-Riemann surface of $K$ over $A$}.

 When no confusion can arise, we  simply denote  by $B_F$ the open set $B_F^{{ Z}}$, and by $B_x$  the open set $B_{\{x\}}$,  for $x \in K$.

  \medskip
  
  Let $R$ be a commutative ring and let { $X:=\sss(R)$} denote  the 
collection 
of prime ideals of $R$.  On { $X$,} we can consider
the {\it Zariski  topology} by taking as closed sets the collection of all sets $V(I):= \{P \in \sss(R) \mid  I
\subseteq P \}$ where $I$ is an ideal of $R$. We denote by  {$X^{\mbox{\tiny{\texttt{zar}}}}$} the prime spectrum of $R$ endowed with the Zariski topology.  If we set $D_a:=\sss(R)\w V(a)$ for all $a \in R$,   it is well known that the family $\{ D_a \mid a \in R\}$ is a basis for the open sets of $X^{\mbox{\tiny{\texttt{zar}}}}$. 
Zariski's topology has several attractive properties related to the 
geometric aspects of the study of the set of prime ideals 
\cite[Chapter I]{Eisenbud}.  For example, { $X^{\mbox{\tiny{\texttt{zar}}}}$} is 
always quasi-compact. { But,} this topology is very coarse.  For 
example, 
{ $X^{\mbox{\tiny{\texttt{zar}}}}$} is always Kolmogoroff, but almost never 
Hausdorff {(more precisely,  $X^{\mbox{\tiny{\texttt{zar}}}}$ is Hausdorff if and only if $\dim(R) =0$ { \cite[Th\'eor\`eme 1.3]{maroscia-0} or \cite[Theorem 3.6]{gilmer-0}}).}  

\smallskip

Many authors have considered a finer topology on the prime spectrum of a ring, known as 
the 
{\it constructible topology} (\cite [pages 337-339]{EGA} or  
\cite[Chapter 3,  Exercises 27, 28 and 30]{AM})  or as the {\it patch to\-po\-logy}  \cite{ho}. In order to introduce  such a  topology in a more general setting, with a simple set theoretical approach, we need some notation and terminology. Given a topological space $\mathcal X$, { with the notation used in \cite[Section 2]{sch-tr}} we set:
$$
\begin{array}{rl}
\mathring{\mathcal K}(\mathcal X):=&\hskip -5pt  \{ U \mid U \subseteq \mathcal X,\, U \mbox{ open and quasi-compact in } \mathcal X\},\\
\overline{\mathcal K}(\mathcal X):=&\hskip -5pt  \{ \mathcal X \setminus U \mid U\in \mathring{\mathcal K}(\mathcal X) \},\\
{\mathcal K}(\mathcal X):=&\hskip -5pt  \mbox{the Boolean algebra of the subsets of $\mathcal X$ generated by $\mathring{\mathcal K}(\mathcal X)$}, \\
\end{array}
$$
 i.e.,  ${\mathcal K}(\mathcal X)$ is the smallest subset of $\mathscr{B}(\mathcal X)$ containing $\mathring{\mathcal K}(\mathcal X)$ and closed with respect to $\cup$\ \!, $\cap$\ \!, and complementation.   As in \cite{sch-tr}, we call the {\it constructible topology on $\mathcal X$} the topology on $\mathcal X$ having ${\mathcal K}(\mathcal X)$ as a basis {(for the open sets). 
We denote by ${\mathcal X}^{\mbox{\tiny{\texttt{cons}}}}$ the set $\mathcal X$ equipped with the constructible topology and we call {\it constructible sets of $\mathcal X$} the elements of ${\mathcal K}(\mathcal X)$  (for Noetherian topological spaces, this notion coincides with that given in \cite[\S 4]{ch}) and {\it proconstructible sets} the closed sets of ${\mathcal X}^{\mbox{\tiny{\texttt{cons}}}}$.

If $X :=\sss(R)$ for some ring $R$, then it is well known that the Zariski topology on $X$ has the set $ \mathring{\mathcal K}(X)$ as a basis (for the open sets) and thus the constructible topology on $X$ is a refinement of the Zariski  topology.  On the other hand, the \emph{constructible topology on $\sss(R)$} is the topology having  the constructible subsets as subsets that are simultaneously  open and closed \cite[({\bf I}.7.2.11) and ({\bf I}.7.2.12)]{EGA}.
More precisely,

\begin{prop} \label{sharp-spec} \cite[Chapter 3, Exercise 28]{AM}. 
Let $R$ be a ring and $X :=\sss(R)$.
Denote by $X^{\mbox{\footnotesize{\texttt{\#}}}}$ the set $X$ endowed with {\rm  the {\texttt{\#}}-topology}, defined as the coarsest  topology on $X$ in which the subsets of type $D_a$ where $a\in R$
are both open and closed.
Then, \  $X^{\mbox{\tiny{\texttt{cons}}}} =X^{\mbox{\footnotesize{\texttt{\#}}}}$. 
In particular, the constructible
topology on $X$ is 
Hausdorff and so $X^{\mbox{\tiny{\texttt{cons}}}} $ is a compact space.   
\end{prop}

\begin{oss}
{\bf (a)} It can be easily shown that the constructible topology on the prime spectrum of a ring is the coarsest topology having as closed sets the closed sets and the quasi-compact open sets  of the Zariski topology, i.e., the family of sets $\{V(I),\ D_a  \mid I \mbox{ is an ideal of } R,\  a \in R\}$ is a subbasis for the closed subspaces of  $\sss(R)^{\mbox{\tiny{\texttt{cons}}}}$  \cite[\S 2, page 45]{ho}.  

{\bf (b)} Another way to describe the constructible topology on $\sss(R)$ is given by taking as closed sets 
the collection of all subsets of $\sss(R)$ of the form $\{f^{-1}(Q) \mid Q \in \sss(S)\}$, where  $f: R \rightarrow S$ is any ring homomorphism \cite[Exercise 27, page 48]{AM}.

{\bf (c)} Let $X=\sss(R)$. If $X^{\mbox{\tiny{\texttt{zar}}}}$ is a Noetherian spectral space, the constructible sets  of $X$ are exactly the finite unions of locally closed subspaces (i.e.,  subspaces obtained by intersection of a closed set with an open set of $X^{\mbox{\tiny{\texttt{zar}}}}$)  
\cite[({\bf 0}.2.3.11) and ({\bf 0}.2.4.1)]{EGA}. By a well known result by Chevalley, if $f: R \rightarrow T$ is a ring homomorphism of finite type and $R$ is a Noetherian ring, then $\{f^{-1}(Q) \mid Q \in \sss(T)\}$ is a constructible subset of $\sss(R)$ \cite[Corollary 14.7]{Eisenbud}.

{\bf (d)} Note also that if $Y$ is a subset of  $X:=\sss(R)$ and if ${Y}^{^{\!\uparrow}} := \{ P \in  X\mid P \supseteq Q, \mbox{ for some } Q \in Y \}$, then the closure of a subset $Y$  of $X$ in the Zariski topology and in the constructible topology are related by the following formula \cite[Lemma (1.1)]{fo}:
$$ \ad^{\mbox{\tiny{\texttt{zar}}}}(Y) = (\ad^{\mbox{\tiny{\texttt{cons}}}}(Y))^{^{\!\uparrow}}.$$

\end{oss}

Recently, Fontana and Loper in \cite{folo}  have considered ``another'' topology on $X:=\sss(R)$ by 
using the notion of an ultrafilter.
Let $C$ be a subset of $X$, and 
let $\uu$ be an ultrafilter on the set $C$.   Set  
 $$
P_{\uu} 
:= \{a \in R  \mid V(a) \cap C \in \uu \}.
$$
  By an argument similar 
to that used in \cite[Lemma 2.4]{calota}, it can be easily shown that 
$P_{\uu}$ is a prime ideal of $R$.     We 
call $P_{\uu}$ an {\it ultrafilter limit point of $C$ in $X$}. This 
notion of 
ultrafilter limit points of collections of prime ideals has been used to 
great effect in 
several recent papers  \cite{calota}, \cite{Loper1}, and  
\cite{Loper2}.  If $\uu$ is a trivial ultrafilter on $C$ then, by definition, there is a 
prime $P \in C$ such that { $\uu = \{Z \in\ms B(C) \mid P \in Z\} \ ( =: \beta^P_C)$ and it  is straightforward in this case that $P_{\uu} = P \in C$ } { \cite[page 2918]{folo}.}  On the 
other hand, if $\uu$ is nontrivial, then it is not at all 
clear that the prime ideal $P_{\uu}$ should lie in $C$.   That 
motivates the following definition.
 Let $R$, $X$ and $C$ be as above.  We say that the set $C$ is 
\it ultrafilter closed \rm  in $X$ if it contains all of its ultrafilter 
limit  points. 
It is not hard to see that the ultrafilter closed subsets 
of $X$ define a topology on the set $X$, called {\it the 
ultrafilter topology on $X$} {\cite[Definition 1]{folo}.}  We denote by $X^{\mbox{\tiny{\texttt{ultra}}}}$  the set  of prime ideals of $R$ endowed with the ultrafilter topology.
One of the main results of a recent paper by Fontana and Loper is the following.

\begin{thm}\label{folo} {\rm \cite[Theorem 8]{folo}}. 
 Let $R$ be a commutative ring and let $X:=\sss(R)$.  Then,  $X^{\mbox{\tiny{\texttt{ultra}}}} = X^{\mbox{\tiny{\texttt{cons}}}}$  (i.e., the ultrafilter topology coincides with the constructible topo\-lo\-gy on the prime spectrum of a ring).
\end{thm}

\section{The ultrafilter topology on $\zar(K|A)$}
\label{const}

 Let $K$ be a field and $A$ a subring of $K$.  Taking as starting point the situation on the prime spectrum of a ring, the next goal is a study of  some topologies on the space ${ Z:= } \zar(K|A)$ that are finer than the Zariski topology.
 
 We start by recalling a very useful fact.

\begin{prop}\label{caeloptar}
Let $K$ be a field and $A$ a subring of $K$. If $Y$ is a nonempty subset of $ Z:=  \zar(K|A)$ and $\ms U$ is an ultrafilter on $Y$, then $A_{\ms U, Y}:=A_{\ms U}:=\{x\in K \mid B_x\cap Y\in \ms U\}$ { is a valuation domain belonging to $Z$.}
\end{prop}
\noindent
\textsc{Proof}. By \cite[Lemma (2.9)]{calota}, $A_{\ms U}$ is a valuation ring of $K$. It remains to show that $A\subseteq A_{\ms U}$. This follow immediately noting that, for every $x\in A$, we have $B_x=\zar(K|A)$, and hence $B_x\cap Y=Y\in \ms U$.\hfill$\Box$

\begin{oss}\label{trivial} The previous statement shows that, 
 if $Y \subseteq Z := \zar(K|A)$, we have a canonical map:
$$
\pi_Y: \boldsymbol{\beta}(Y) \rightarrow  Z\,, \quad \ms U \mapsto A_{\ms U, Y}:=
\{x\in K \mid B_x\cap Y\in \ms U\},
$$ 
and, in this case, $Y \subseteq \im(\pi_Y)$, since  for each $V \in Y$, taking the trivial ultrafilter ${\beta_Y^V}  \in \boldsymbol{\beta}(Y)$, { we have $A_{\beta_Y^V, Y} =V$}.
\end{oss}

The previous remark leads naturally to the following crucial definition of this section.
Let $K$ be a field and $A$ a subring of $K$. A subset $Y$ of   $\zar(K|A)$  is called \emph{stable for ultrafilters} if,  for each $\ms U\in  \boldsymbol{\beta}(Y) $, $A_{\ms U\!, Y}\in Y$ (or, equivalently, with the notation of Remark \ref{trivial}, $\im(\pi_Y) =Y$).

\begin{prop}\label{topologia}
Let $K$ be a field, $A$ a subring of $K$ and $Z := \zar(K|A)$. Then, the collection of all subsets of { $Z$} stable for ultrafilters is the family of closed sets for a topology on $Z$ called {\rm the ultrafilter topology of the Zariski-Riemann surface} {$Z$}.\end{prop}
\noindent
\textsc{Proof}. The empty set and { $Z$} are clearly stable for ultrafilters. Now, consider two subsets $C',C''$ of {$Z$} stable for ultrafilters, set $Y:=C'\cup C''$, and let $\ms U$ be an ultrafilter on $Y$. By Lemma \ref{siparte}(5), we can assume, without loss of generality, that $C'\in \ms U$. Then $\ms U' := \ms U^{{C'}}:=\{Z\cap C' \mid Z\in \ms U\}$ is an ultrafilter on $C'$, by Lemma \ref{siparte}(3). We want to show that $A_{\ms U}=A_{\ms U'}$. { Let } $x\in A_{\ms U'}$. 
Then $B_x\cap C'\in \ms U' \subseteq \ms U$ (by Lemma \ref{siparte}(3)). 
Since $B_x\cap C'\subseteq B_x\cap Y$, it follows immediately that $B_x\cap Y\in \ms U$ and hence $x\in A_{\ms U} \ (=\{x\in K \mid B_x\cap Y\in \ms U\})$. 
Conversely, { let } $x\in A_{\ms U}$. Since $B_x\cap Y\in \ms U$, we have $B_x\cap C'=(B_x\cap Y)\cap C'\in \ms U'$. Hence, $x\in A_{\ms U'} \ (= \{x\in K \mid B_x\cap C' \in \ms U'\})$ { and so } $A_{\ms U}=A_{\ms U'}$.  As $C'$ is stable for ultrafilters, we have $A_{\ms U}= A_{\ms U'}\in C'\subseteq Y$ and so $Y$ is also stable for ultrafilters. By induction, we { easily deduce } that the union of a finite family of subsets stable for ultrafilters  is still stable for ultrafilters. 
Now, let $\ms C$ be any collection of subsets stable for ultrafilters in $Z$ and set $Y:=\bigcap \ms C$. { Let  $\ms U$ be an ultrafilter on $Y$.} 
For every $C\in \ms C$, clearly $Y \subseteq C$ and so, by Lemma \ref{siparte}(4), $\ms U{_{\!{\!{C}}}}:=\{W\in \ms B(C) \mid W\cap Y\in \ms U\}$ is an ultrafilter on $C$.
 Moreover, as before, it is easily seen that $A_{\ms U}=A_{{\ms U}{_{\!{\!{C}}}}} \in C$. This proves that $A_{\ms U}\in \bigcap \ms C$, and thus  every intersection of subsets { of $Z$ } stable for ultrafilters is still stable for ultrafilters. \hfill$\Box$
 
 \medskip

As above, let $Z:= \zar(K|A)$, we denote by $Z^{\mbox{\tiny{\texttt{ultra}}}}$ [respectively, $Z^{\mbox{\tiny{\texttt{cons}}}}$; $Z^{\mbox{\tiny{\texttt{zar}}}}$] the space of valuation domains of $K$ containing $A$ equipped with the ultrafilter  topology [respectively, with the constructible topology; with the Zariski topology].  The next goal is to  compare $Z^{\mbox{\tiny{\texttt{ultra}}}}$ with $Z^{\mbox{\tiny{\texttt{cons}}}}$ and  $Z^{\mbox{\tiny{\texttt{zar}}}}$.

\begin{thm}\label{fine}
Let $K$ be a field, $A$ a subring of $K$ 

and let $Z:= \zar(K|A)$.
\begin{enumerate}
\item[(1)]   The ultrafilter topology   is finer than the Zariski topology on $Z$.
\item[(2)] For any  subset $S$ of $K$, $B^Z_S \ (:= B_S := \{V \in Z \mid V \supseteq S\})$ is a closed set in  the ultrafilter topology. 
In particular, the basic open  sets of the Zariski topology of $Z$ are both open and closed in the ultrafilter topology.
\item[(3)]  We denote by $Z^{\mbox{\footnotesize{\texttt{\#}}}}$ the set $Z$ endowed with {\rm the \texttt{\#}-topology}, defined as the coarsest topology  for which the set $B_F$ is both open and closed, for every finite subset $F$ of $K$.  Then, $Z^{\mbox{\footnotesize{\texttt{\#}}}}$ is a Hausdorff topological space. 
\item[(4)]  The \texttt{\#}-topology on $Z$  is the coarsest topology having as closed sets the closed sets and the quasi-compact open sets  of $Z^{\mbox{\tiny{\texttt{zar}}}}$, i.e.,  

\begin{center}
$
{\boldsymbol{\mathcal C}}^{\mbox{\footnotesize{\texttt{\#}}}} := \left\{ B_F\ ; \, \bigcap \{Z\w B_G \mid G \in \mathcal G\} \mid  F \in \ms B_{\mbox{\Tiny{{\rm\texttt{fin}}}}}(K) , \,  \mathcal G\subseteq \ms B_{\mbox{\Tiny{\rm\texttt{fin}}}}(K) \right\}
$
\end{center}
 is a subbasis for the closed  subsets of  $Z^{\mbox{\footnotesize{\texttt{\#}}}}$.   
 \item[(5)]  $Z^{\mbox{\tiny{\texttt{ultra}}}}$  is a (Hausdorff) compact topological space.
\item[(6)]  $Z^{\mbox{\tiny{\texttt{ultra}}}} = Z^{\mbox{\footnotesize{\texttt{\#}}}}  = Z^{\mbox{\tiny{\texttt{cons}}}}$.
\end{enumerate} 
\end{thm}
\noindent
\textsc{Proof}. 
(1)  Since $\{B_F\mid F\in \ms B_{\mbox{\Tiny{\texttt{fin}}}}(K)\}$ is a { basis for the open sets on $Z^{\mbox{\tiny{\texttt{zar}}}}$,} it is enough to prove that $Z \w B_F$ is stable for ultrafilters, for every $F\in \ms B_{\mbox{\Tiny{\texttt{fin}}}}(K)$.
{  Assume,} by contradiction, that there exists an ultrafilter $\ms U$ on $Y:=Z \w B_F$ such that $A_\ms U\notin Y$. It follows that $F\subseteq A_{\ms U}$, and then $B_x\cap Y\in \ms U$, for every $x\in F$. Since $F$ is finite, we have $B_F\cap Y\in \ms U$. This is a contradiction by the definition of $Y$ { (and by the fact that $\emptyset$ does not belong to any filter).} 

 (2)  Apply Proposition \ref{caeloptar}, after observing   that $B_S=\zar(K|A[S])$.

(3) Let $V$ and $W$ be two distinct elements of $Z$, and, without loss of generality, we can take an element $x\in V\w W$. By assumption, the sets $B_x$ and $Z\w B_x$ are disjoint open neighborhoods of $V$ and $W$, respectively, in the topological space $Z^{\mbox{\footnotesize{\texttt{\#}}}}$.

(4) It is clear that each set in ${\boldsymbol{\mathcal C}}^{\mbox{\footnotesize{\texttt{\#}}}}$ is closed in the \texttt{\#}-topology and every topology in which the sets of type 
 $B_F$ (for $F \in \ms B_{\mbox{\Tiny{\texttt{fin}}}}(K)$) are both open and closed  must be finer than the topology having ${\boldsymbol{\mathcal C}}^{\mbox{\footnotesize{\texttt{\#}}}}$ as subbasis for the closed sets. Conversely, it is obvious that, in this last topology, each set of type 
 $B_F$ (for $F \in \ms B_{\mbox{\Tiny{\texttt{fin}}}}(K)$) is both open and closed.

 (5) First, we note that,  by (2) and (3), the ultrafilter topology on $Z$ is finer than the \texttt{\#}-topology and so  $Z^{\mbox{\tiny{\texttt{ultra}}}}$ is a Hausdorff space. 
Let $\ms C$ be a collection of closed subsets of $Z^{\mbox{\tiny{\texttt{ultra}}}}$ with the finite intersection property. By Lemma \ref{siparte}(1 and 2), we can { find} an ultrafilter $\ms U$ on $Z$ containing $\ms C$. Now, take a closed set $C\in \ms C$, and consider the ultrafilter $\ms U^{{\!}C}\in \boldsymbol{\beta}(C)$ induced by $\ms U$ (Lemma \ref{siparte}(3)).
 By the same argument used in the proof of Proposition \ref{topologia}, we have $A_{\ms U}=A_{\ms U^{{\!}C}}$. Keeping in mind that every element of $\ms C$ is stable for ultrafilters, we deduce that  $A_{\ms U}\in \bigcap \ms C$, and so $\bigcap \ms C \neq \emptyset$.
 
(6)  By (2), the identity map ${\texttt{id}}_Z: Z^{\mbox{\tiny{\texttt{ultra}}}} \rightarrow Z^{\mbox{\footnotesize{\texttt{\#}}}}$ is continuous.
 Moreover, since $Z^{\mbox{\tiny{\texttt{ultra}}}} $ is compact (by { (4)}) 
 and $Z^{\mbox{\footnotesize{\texttt{\#}}}}$ is Hausdorff (by { (3)}), \texttt{id}$_Z$ is a closed map (cf., for instance, \cite[Chapter IX, Theorem 2.1]{Du}), and hence  is a homeomorphism (cf., for instance, \cite[Chapter III, Theorem 12.2]{Du}).  { Finally, the equality 
$Z^{\mbox{\footnotesize{\texttt{\#}}}}=Z^{\mbox{\tiny{\texttt{cons}}}} $
 follows immediately from (4) and from the definition of the constructible topology}. 
  \hfill$\Box$

\begin{oss} {\rm Note that, {\sl mutatis mutandis},  the proofs of points (5) and (6) of the previous theorem provide another very short and purely topological proof  of the fact that the ultrafilter topology and the patch (or, constructible) topology coincide on the prime spectrum of a ring \cite[Theorem 8]{folo}. The idea for this type of topological argument  was already in \cite[Appendix, Theorem 3.12]{Finocchiaro}. A similar (topological) proof, in the case of the prime spectrum of a ring, was given independently in  \cite{ruvi}. 
}
\end{oss}

\smallskip
From Theorem \ref{fine} ((1) and (5)) (and straightforward topological arguments),  we  easily reobtain  the following well-known fact  \cite[Chapter VI, Theorem (40)]{zasa}.

 \begin{cor} \label{zariski}
  Let $K$ be a field, $A$ a subring of $K$ { and $Z:=\zar(K|A)$. 
Then,
$Z^{\mbox{\tiny{\texttt{zar}}}}$} is a Kolmogoroff  quasi-compact topological space.
 \end{cor}

\begin{prop}  \label{ad}
Let $K$ be a field,  $A$ a subring of $K$ { and $Z:=\zar(K|A)$. Denote by $
\ad^{\mbox{\tiny{\texttt{ultra}}}}(Y)$ the closure of a subset $Y$  in $Z^{\mbox{\tiny{\texttt{ultra}}}}$.}  Then, 
$
\ad^{\mbox{\tiny{\texttt{ultra}}}}(Y)=\left\{A_{\ms U}\mid \ms U\in {\boldsymbol{\beta}}(Y)\right\}.
$
\end{prop}
\noindent
\textsc{Proof}. We begin noting that, by  { Theorem \ref{fine}((4)  and (6)),} a basis for the open sets { of } the ultrafilter topology is given by
$$
\mathcal B^{\mbox{\tiny{\texttt{ultra}}}}:=\mathcal B:= 
\left\{
B_F\,;\; B_F\cap \left(\bigcap_{i=1}^n(Z\w  B_{F_i})\right) \mid 
F,F_1,F_2, \z, F_n\in \ms B_{\mbox{\Tiny{\texttt{fin}}}}(K),\ n\geq 1 
\right\}.
$$
Now, let $\ms U$ be an ultrafilter on $Y$ and $U$ be an open neighborhood of $A_{\ms U}$ in ${Z}^{\mbox{\tiny{\texttt{ultra}}}} $. By the above remark, we can assume, without loss of generality, that $U$ is of the form $B_F$ or $B_F\cap \bigcap_{i=1}^n(Z\w  B_{F_i})$, for some collection of  finite subsets $F, F_1, F_2, \z,F_n$ of $K$ and some $n\geq 1$. 
If $U=B_F$, then $F\subseteq A_{\ms U}$, 
and so $B_F\cap Y\in \ms U$, by the definition of $A_{\ms U}$.
 In particular, $B_F\cap Y\neq \emptyset$.
  If $U=B_F\cap \bigcap_{i=1}^n(Z\w  B_{F_i})$ we have $B_F\cap Y\in \ms U$,
   by the same argument given above. 
   Moreover, it can be easily shown that $ B_{F_i} \cap Y \notin \ms U$, for each $i$, and hence $ \bigcap_{i=1}^n (Z\w B_{F_i}) \cap Y \in \ms U$. 
     Since $\emptyset $ does not belong to any ultrafilter, it follows that $U\cap Y\neq \emptyset$. This proves that $\{A_{\ms U}\mid \ms U\in {\boldsymbol{\beta}}(Y)\} \subseteq \ad^{\mbox{\tiny{\texttt{ultra}}}}(Y)$. 
Conversely, let $V$ be a valuation domain in $\ad^{\mbox{\tiny{\texttt{ultra}}}}(Y)$. 
If $F$ is a finite subset of $V$ and $F_1,F_2, \z,F_n$ are finite subsets of $K$ such that $F_i\nsubseteq V$, for $i=1, 2, \z, n$,  $B_F\cap \bigcap_{i=1}^n(Z\w  B_{F_i})\cap Y$ is nonempty. Then, it follows immediately that the following family of sets
$$
\begin{array}{rl}
\mathcal B_V:= {\boldsymbol\{}B_F\cap Y, & \hskip -5pt B_F\cap \bigcap_{i=1}^n(Z\w  B_{F_i})\cap Y \mid  \\
& \hskip 40pt F\in \ms B_{\mbox{\Tiny{\texttt{fin}}}}(V),\ F_1, F_2, \z, F_n\in \ms B_{\mbox{\Tiny{\texttt{fin}}}}(K)\w \ms B(V), \,  n\geq 1 {\boldsymbol \}}
\end{array}
$$
is a collection of subsets of $Y$ with the finite intersection property, and thus there exists an ultrafilter $\ms U\in \boldsymbol{\beta}(Y)$ such that $\mathcal B_V\subseteq \ms U$ (Lemma \ref{siparte} (1 and 2)). 
It is enough to show that $A_{\ms U}=V$. If $x\in A_{\ms U} \w V$, then we have 
$(Z\w B_x) \cap Y\in \mathcal B_V\subseteq \ms U$, by construction and, moreover, 
$ B_x \cap Y \in \ms U$, by the definition of $A_{\ms U}$, { which is a contradiction.} 
Conversely, let $x\in V$. Then $V\in B_x$ and, thus, $B_x\cap Y\in \mathcal B_V\subseteq \ms U$. In other words, $x\in A_{\ms U}$. This proves the statement. \hfill$\Box$

\medskip

As  is well known, if $K$ is a field and $A$ is a subring of $K$, we can construct a map $\gamma:\zar(K|A)\rightarrow \sss(A)$ sending a valutation ring $V\in \zar(K|A)$, with maximal ideal  $M_V$, to the prime ideal  $M_V\cap A$ of $A$,  called {\it the center of $V$ over $A$}. { It is well known (by an application of Zorn's Lemma) that $\gamma$ is   a surjective map. }

Moreover, if we consider $Z:=\zar(K|A)$ and $X:=\sss(A)$  as topological spaces both endowed with the Zariski topology  
then, by \cite[Lemma (2.1)]{dofefo},  the map 
$\gamma:Z^{\mbox{\tiny{\texttt{zar}}}} \rightarrow X^{\mbox{\tiny{\texttt{zar}}}}$
  is continuous, since {$\gamma^{-1}(D_a) = B_{a^{-1}}$},
   for each nonzero { $a \in A$.}
   { Moreover,  $\gamma: Z^{\mbox{\tiny{\texttt{zar}}}} \rightarrow X^{\mbox{\tiny{\texttt{zar}}}}$ 
   is also a closed  map, essentially by \cite[Theorem (2.5)]{dofefo} { see also Remark \ref{singleton}). } 
   In particular, $\gamma: Z^{\mbox{\tiny{\texttt{zar}}}} \rightarrow X^{\mbox{\tiny{\texttt{zar}}}}$ is a homeomorphism 
   if and only if $\gamma$ is injective (i.e., if and only if for each 
   $P \in \sss(A)$ there exists a unique valuation domain of $K$ dominating $A_P$).}
    In particular, if $A$ is a Pr\"ufer domain with quotient field $K$, then  $\gamma: Z^{\mbox{\tiny{\texttt{zar}}}} \rightarrow X^{\mbox{\tiny{\texttt{zar}}}}$ 
    is a homeomorphism.

    \begin{oss}\label{singleton} 
\rm Note that, in \cite{dofefo}, the authors consider the case where $A$ is an integral domain with quotient field $K$.
 If $A$ is a subring, but not a subfield of $K$, and if the quotient field of $A$ is a proper subfield of $K$, then we can take the integral closure $\bar{A}$ of $A$ in $K$. In this situation, $\bar{A}$ is an integral domain such that $\zar(K|A)= \zar(K|\bar{A})$.    If $A$ is a subfield of $K$, then $\sss(A)$ is a (discrete) topological space consisting of just one point and so, trivially,  the map $\gamma:Z^{\mbox{\tiny{\texttt{zar}}}} \rightarrow X^{\mbox{\tiny{\texttt{zar}}}}$ in this case is continuous, surjective and closed. 
    \end{oss}


\medskip

The next goal is to study the map $\gamma$ when $Z:=\zar(K|A)$  and $X:=\sss(A)$ are both equipped with the ultrafilter topology (or, equivalently, with the constructible topology (Theorem \ref{folo})).

\begin{thm}\label{continuous-closed}
Let $K$ be a field and $A$ a subring of $K$. Then, the surjective map $\gamma: \zar(K|A)^{\mbox{\tiny{\texttt{ultra}}}} \rightarrow \sss(A)^{\mbox{\tiny{\texttt{ultra}}}}$  is continuous and closed.
\end{thm}
\noindent
\textsc{Proof}. Set as usual $Z:=\zar(K|A)$ and $X:= \sss(A)$. Since $Z^{\mbox{\tiny{\texttt{ultra}}}}$ is compact, by { Theorem \ref{fine}(5)}, and $X^{\mbox{\tiny{\texttt{ultra}}}}= X^{\mbox{\tiny{\texttt{cons}}}}$ is Hausdorff (and compact), by standard topological properties (cf., for instance, \cite[Chapter XI, Theorem 2.1]{Du}), it is enough to show that $\gamma$ is continuous.
 Let $C$ be a closed subset of $X^{\mbox{\tiny{\texttt{ultra}}}}$, $\ms U$ an ultrafilter on $\gamma^{-1}(C) \ (\subseteq Z) $  and let  $\delta:\gamma^{-1}(C)\longrightarrow C$ be the restriction of $\gamma$ to $\gamma^{-1}(C)$. 
 By Lemma \ref{siparte}(4), the collection of sets $$\ms V:=\ms U_{_{\!\delta}} :=\{V\subseteq C\mid \delta^{-1}(V)\in \ms U\}=\{V\subseteq C \mid \gamma^{-1}(V)\in \ms U\}$$ is an ultrafilter on $C$. 
 So, we can consider  ${A_{\ms U}} \in  Z$ (more precisely, ${A_{\ms U}}  =\{x\in K \mid B_x\cap \gamma^{-1}(C)\in \ms U\}$  is a point in the closure of $\gamma^{-1}(C)$ in $Z^{\mbox{\tiny{\texttt{ultra}}}}$, by Proposition \ref{ad}) and { we can also consider the ultrafilter limit point 
 $$
 P_{\ms V}:=\{a \in A \mid V(a)\cap C\in  \ms V\}
 $$ 
 which is a prime ideal of $A$ \cite[page 2918]{folo}.}
  We claim that  the center on $A$ (of the maximal ideal $M{_{\!\ms U}}$) of the valuation domain ${A_{\ms U}}$ coincides with $P_{\ms V}$, i.e.,  $M{_{\!\ms U}}\cap A=P_{\ms V}$. 
  As a matter of fact, let $a\in P_{\ms V}$. By definition, it follows immediately that $ \gamma^{-1}(V(a)\cap C) = \gamma^{-1}(V(a))\cap \gamma^{-1}(C)\in \ms U$.
   Now assume, by contradiction, that $a$ is a unit in $A_{\ms U}$. 
   Equivalently, $B_{a^{-1}}\cap \gamma^{-1}(C)$ belongs to $\ms U$. Since $\gamma^{-1}(V(a))\cap B_{a^{-1}}\cap \gamma^{-1}(C) \in \ms U$, in particular, $\gamma^{-1}(V(a))\cap B_{a^{-1}}\cap \gamma^{-1}(C)$ is nonempty. 
   Therefore,  there exists a valuation domain $W\in \gamma^{-1}(C)$ such that  $a^{-1}\in W$ and $a\in \gamma(W) := M_W\cap A$, where $M_W$ is the maximal ideal of the valuation domain $W$.
  It follows immediately  that $1\in M_W$, a contradiction. Therefore, $ P_{\ms V} \subseteq M{_{\!\ms U}} \cap A$.
  Conversely, let $a\in M_{A_{\ms U}}\cap A$, $a \neq 0$.
    Then, in particular,  $a^{-1}\notin A_{\ms U}$ and, since as we have already observed $\gamma^{-1}(D_a) = B_{a^{-1}}$, we have 
     $B_{a^{-1}}\cap \gamma^{-1}(C)=\gamma^{-1}(D_a\cap C)\notin \ms U$. Hence $D_a\cap C\notin \ms V$  and, thus,  finally $V(a)\cap C\in \ms V$, since $\ms V$ is an ultrafilter on $C$ (Lemma \ref{siparte}(5)). Therefore $a \in  P_{\ms V}$.
    This shows that  $M{_{\!\ms U}}\cap A=P_{\ms V}$. 
    
    Since,  by \cite[Theorem 8]{folo}, $C$ is stable for ultrafilters,  we have $\gamma (A_{\ms U})=P_{\ms V}\in C$,  and so $A_{\ms U} \in  \gamma^{-1}(C)$. Therefore, we deduce that $\gamma^{-1}(C)$ is closed in $Z^{\mbox{\tiny{\texttt{ultra}}}}$ and so
    the conclusion follows.\hfill$\Box$

      \begin{oss} \label{gamma-prufer} \rm
     Note that, with the notation of the previous Theorem \ref{continuous-closed} (and its proof), if $A$ is a Pr\"ufer domain, the map $\gamma: Z^{\mbox{\tiny{\texttt{ultra}}}} \rightarrow X^{\mbox{\tiny{\texttt{ultra}}}}$ is a homeomorphism, since in this case (as observed just before Remark  \ref{singleton}) $\gamma$ is injective.
      \end{oss}

\begin{center}
ACKNOWLEDGEMENTS
\end{center}

The authors are grateful to the referee for his/her useful comments and for pointing out the publication of a recent related paper by Luz M. Ruza and Jorge Vielm \cite{ruvi}.

 

\begin{thebibliography}{9999}
\footnotesize

\bibitem{AP}
Victor Alexandru and Nicolae Popescu,  Sur une classe de prolongements \`a $K(X)$ d'une valuation sur un corps $K$, {\it Rev. Roumaine Math. Pures Appl.} {\bf 33} (1988), 393--400.


\bibitem{APZ}
Victor Alexandru, Nicolae Popescu and Alexandru Zaharescu,  All valuations on $K(X)$, {\it J. Math. Kyoto Univ.} {\bf 30} (1990),  281--296.




\bibitem{AM}  Michael F. Atiyah and Ian G. Macdonald, {\it Introduction to 
commutative algebra,} Addison-Wesley, Reading MA, 1969.


\bibitem{calota}
Paul-Jean Cahen, Alan Loper, and Francesca Tartarone, Integer-valued polynomials and Pr\"ufer $v-$multiplication domains, {\it J. Algebra} {\bf 226} (2000), 765--787.

\bibitem{ch}
Claude Chevalley et Henri Cartan, {Sch\'emas normaux; morphismes; ensembles constructibles,} {\it S\'eminaire Henri Cartan}  {\bf 8} (1955-1956), Exp. No. 7, 1--10.

\bibitem{dofefo} David E. Dobbs, Richard Fedder, and Marco Fontana, Abstract Riemann surfaces of integral domains and spectral spaces. {\it Ann. Mat. Pura Appl.}   {\bf 148 } (1987), 101--115.

\bibitem{dofo}
David E. Dobbs and  Marco Fontana, Kronecker Function Rings and Abstract Riemann Surfaces, { \it J.
Algebra} {\bf 99} (1986), 263--274.




\bibitem{Du}
James Dugundji, {\it Topology},  Allyn and Bacon, Boston, 1966.


\bibitem {Eisenbud}
David Eisenbud, {\it Commutative algebra with a view toward algebraic 
geometry}, Springer, Berlin, 1994.


\bibitem{Finocchiaro}
Carmelo Finocchiaro, {\it Amalgamation of algebras and the ultrafilter topology on the Zariski space of valuation overrings of an integral domain'}, Ph.D. Thesis, Universit\`a degli Studi ``Roma Tre'', December 2010.


\bibitem{fi-fo-lo-11}
Carmelo Finocchiaro, Marco Fontana and K. Alan Loper, {\it The ultrafilter topology on spaces of valuation domains}. Submitted.


\bibitem{fo}
Marco  Fontana, Topologically defined classes of commutative rings,  {\it Ann. Mat. Pura Appl.} {\bf 123} (1980), 331--35.




\bibitem{folo} Marco Fontana and  K. Alan Loper, The patch topology and the ultrafilter topology on the prime spectrum of a commutative ring, {\it Comm. Algebra} {\bf 36} (2008), 2917--2922.



\bibitem{gilmer-0}
R. Gilmer,   Background and preliminaries on zero-dimensional rings, in {\sl ``Zero-dimensional Commutative Rings''}, David F. Anderson (Editor), David Dobbs (Editor), M. Dekker Lecture Notes in Pure and Applied Mathematics, {\bf 171}, 1995.


\bibitem{EGA}
Alexander Grothendieck   et   Jean Dieudonn\'e, {\it \'El\'ements de G\'eom\'etrie 
Alg\'ebrique I}, Springer, Berlin, 1970.




\bibitem{ho} Melvin Hochster, Prime ideal structure in commutative rings, {\it Trans. Amer. Math. Soc.} {\bf 142} (1969),   43--60.



\bibitem{hu-kn}
Roland Huber and Manfred Knebusch, On valuation spectra, in
{\sl ``Recent advances in real algebraic geometry and quadratic forms: proceedings of the RAGSQUAD year''}, Berkeley, 1990-1991, {\it Contemp. Math.} {\bf 155}, Amer. Math. Soc. Providence RI, 1994].

\bibitem{jech} Thomas Jech, {\it Set Theory}, Springer, New York, 1997 (First Edition, Academic Press, 1978).



\bibitem{KMPV}
Shigeru Kobayashi, Hidetoshi Marubayashi, Nicolae Popescu and  Constantin Vraciu,  Total valuation rings of $K(X,\sigma)$ containing $K$, {\it Comm. Algebra} {\bf 30} (2002), 5535--5546.


\bibitem{ku}
Franz-Viktor Kuhlmann, Places of algebraic fields in arbitrary characteristic, {\it Advances Math.} {\bf 188} (2004), 399--424.


\bibitem{Loper1} 
K. Alan Loper {Sequence domains and integer-valued polynomials},  
{\it J. Pure  Appl. Algebra}   {\bf 119}  (1997), 185--210.

\bibitem{Loper2}  
K. Alan Loper {A classification of all  $D$ such that  Int($D$) 
is a Pr\"ufer 
domain},  {\it Proc. Amer. Math. Soc.}
 {\bf  126}   (1998),  657--660.
 

 \bibitem{maroscia-0}
 Paolo Maroscia, 
Sur les anneaux de dimension z\'ero. {\it Atti Accad Naz. Lincei Rend. Cl. Sci. Fis. Mat. Natur.} (8) 56 (1974), 451--459.

\bibitem{PV-1}
Nicolae Popescu and Constantin Vraciu,
 On the extension of valuations on a field $K$ to $K(X)$. I. {\it Rend. Sem. Mat. Univ. Padova} {\bf 87} (1992), 151--168.
 
\bibitem{PV-2}
Nicolae Popescu and Constantin Vraciu,  On the extension of a valuation on a field $K$ to $K(X)$. II. {\it Rend. Sem. Mat. Univ. Padova} {\bf 96} (1996), 1--14.

\bibitem{PZ}
Nicolae Popescu and Alexandru Zaharescu,  On a class of valuations on $K(X)$, XI{\texttt th} National Conference of Algebra (Constan\c{t}a, 1994). {\it An. \c{S}tiin\c{t}. Univ. Ovidius Constan\c{t}a,  Ser. Mat.} {\bf 2} (1994), 120--136.
 
 

\bibitem{ruvi}
Luz M. Ruza and Jorge Vielm,
 {\it The equality
of the patch topology and the ultrafilter topology: a shortcut},  Appl. Gen.
Topol. {\bf 12} (2011),  15--16.


 \bibitem{sch-tr}
 Niels Schwartz and Marcus Tressl, Elementary properties of minimal and maximal points in Zariski spectra, {\it J. Algebra} {\bf 323} (2010), 698--728.
 
 
  


 \bibitem{za} 
 Oscar Zariski, {The compactness of the Riemann manifold of an abstract field of algebraic functions},  {\it Bull. Amer. Math. Soc.} {\bf 50} (1944), 683--691.


\bibitem{zasa} Oscar Zariski and Pierre Samuel, {\it Commutative Algebra, Volume 2}, Springer Verlag, Graduate Texts in Mathematics {\bf 29}, New York, 1975 (First Edition, Van Nostrand, Princeton, 1960). 



\end{thebibliography}
\end{document}